\documentclass[a4paper,fleqn]{cas-dc}

\usepackage[numbers]{natbib}
\usepackage{pgfplots}
\usepackage{xcolor}
\usepackage{url}

\usepackage{siunitx}

\def\tsc#1{\csdef{#1}{\textsc{\lowercase{#1}}\xspace}}
\tsc{WGM}
\tsc{QE}
\tsc{EP}
\tsc{PMS}
\tsc{BEC}
\tsc{DE}

\begin{document}
\let\WriteBookmarks\relax
\def\floatpagepagefraction{1}
\def\textpagefraction{.001}
\shorttitle{Batteries and the British Energy System}
\shortauthors{W.Bukhsh}

\title [mode = title]{Batteries and the British Energy System}                      



\author[]{Waqquas Bukhsh}[orcid=0000-0002-5765-0747]
\ead{waqquas.bukhsh@strath.ac.uk}
\ead[url]{https://www.strath.ac.uk/staff/bukhshwaqquasdr/}

\credit{Conceptualization of this study, Methodology, Software}

\affiliation[]{organization={Department of Electronic and Electrical Engineering},
                addressline={University of Strathclyde}, 
                city={Glasgow},
                country={Scotland}}

%
%
%
%



\begin{abstract}
Batteries are becoming a central part of modern energy systems, especially as electricity, transport and heat are decarbonised. In Great Britain, batteries already play an important role by providing flexibility and acting as a buffer for the system, and their importance will continue to grow as the system moves towards net zero by 2050. This perspective reviews how the role of batteries in Great Britain is evolving, their current value, and the growing contribution of batteries in the transport sector. It argues that the full value of electric vehicle batteries will only be realised through coordinated dispatch and better integration into system operation. Evidence from future scenarios developed by the National Energy System Operator shows that while electric vehicles could provide most of the storage capacity, unmanaged charging could create challenges for the system rather than benefits. Unlocking this value will require improvements in market design, optimisation, and control strategies.
\end{abstract}

%

\begin{keywords}
Batteries\sep Energy Storage\sep  Electric Vehicles \sep Ancillary Services \sep Frequency Response \sep Smart Charging of EVs 
\end{keywords}

\maketitle

\section{Introduction}

Energy systems around the world are undergoing their most profound transformation in over a century. Driven by the urgent need to limit climate change, governments worldwide have committed to decarbonisation their energy sectors. The Paris Agreement of 2015 established a global framework for limiting warming to well below $2 \si{\degreeCelsius}$ above pre-industrial levels, and more than 140 countries have since set net zero targets~\cite{unfccc2015paris}. Meeting these commitments requires a complete shift away from fossil fuels and towards renewable energy sources, particularly wind and solar power. Global installed renewable capacity has grown at an unprecedented rate, exceeding 3,000 GW for the first time in 2023~\cite{irena2024renewable}. This shift introduces a fundamental challenge, which is wind and solar generation is variable and weather-dependent. As their share of the electricity mix grows, so too does the need for flexible resources capable of storing surplus energy and releasing it when demand requires. Energy storage has therefore is becoming increasing important~\cite{chu2012opportunities}.

Electricity systems have always needed a way to balance supply and demand. Historically, this role was played predominantly by pumped hydro storage, which remains the largest source of grid-scale storage worldwide, accounting for over 90\% of installed storage capacity globally as recently as 2020~\cite{iea2021storage}. However, pumped hydro is constrained by it's geography and continuously scaling it to meet demand for flexibility is just not possible~\cite{GUERRA2024122447}. Battery storage, and lithium-ion technology in particular, has stepped into this gap. Driven largely by investment from the electric vehicle industry, lithium-ion battery costs fell by more than 90\% between 2010 and 2023~\cite{bnef2023evo}, transforming what was once a niche and expensive technology into a mainstream grid resource~\cite{RoyalSociety2023LargeScale}. Battery energy storage systems are now deployed at scale, with global installed capacity growing at over 80\% per year in recent years~\cite{iea2024battery}.

Great Britain sits at the forefront of this transition. The country has committed to clean power by 2030 and net zero by 2050~\cite{DESNZ2024_clean_power_2030}, targets that require a rapid and sustained expansion of wind and solar generation. Wind and solar already account for a substantial and growing share of GB electricity generation, at times supplying over 70\% of demand. However, managing a system with high levels of renewable generation is costly. It requires more reserve capacity that can respond over different time periods and be ready when needed. There is also a more frequent need for balancing actions. This can arise from forecasting errors or from transmission constraints that limit the ability to move renewable power to where it is needed. As a result, balancing costs in Great Britain have increased steadily. Annual costs rose from £1.2 billion in 2018/19 to £2.7 billion in 2024/25, more than doubling over this period. They reached a peak of £4.15 billion in 2022/23~\cite{NESO_AnnualBalancingReport_2025}. While several factors contribute to this trend, a key issue is the reduction in system flexibility and the reliance on expensive dispatchable resources when alternatives are limited.

This perspective examines the evolving role of battery storage in the Great Britain electricity system, from early grid scale deployments to the large capacities projected in the system operator's NESO Future Energy Scenarios up to 2050~\cite{neso2025fes}. Particular attention is given to electric vehicle batteries, which may represent the largest untapped source of flexibility in the system. While grid connected battery storage is growing rapidly, the combined capacity of millions of electric vehicles could eventually exceed all other sources. However, this potential can only be realised through careful and coordinated management. Without it, uncoordinated charging and discharging could create new risks for system stability. Unlocking this potential is one of the central challenges of the energy transition.

\section{Batteries in Great Britain}

\subsection{Early Development}
Great Britain's first grid scale battery project was delivered by UK Power Networks at Leighton Buzzard in Bedfordshire. The 6 MW/10 MWh system began operation at the end of 2014 under the Smarter Network Storage project~\cite{UKPN_SNS_2016}. The project showed that batteries could provide transmission services through Dynamic Firm Frequency Response. It also confirmed that batteries can respond to frequency changes much faster than conventional generators due to the speed of lithium ion technology.

In 2016 and 2017, the system operator (National Grid ESO at the time) ran a competitive tender for Enhanced Frequency Response, procuring around 200 MW of capacity from battery projects~\cite{nationalgrid2016efr}. These contracts were a turning point, providing the revenue certainty needed to attract investment. This led to a wave of projects between 10 MW and 50 MW across England and Scotland, often located at lower cost grid connection points.

On 9 August 2019, a lightning strike hit a transmission circuit near Eaton Socon, Cambridgeshire, causing the sudden loss of about 1,378 MW of generation. Approximately 475 MW of battery capacity responded immediately, helping to slow the frequency fall~\cite{stoker2019blackout}. Without this response, the demand disconnection affecting over one million customers would likely have been more severe and longer lasting. This event clearly demonstrated the value of battery storage in maintaining system reliability.

\subsection{The Current Fleet}
Building on the early projects and market signals described above, battery storage in Great Britain has expanded rapidly in recent years. What began as a small number of demonstration and pilot projects has now developed into a large and growing fleet of commercial assets.

By the end of 2024, total battery energy storage capacity in Great Britain had reached approximately 4.7 GW. At the same time, the nature of new battery projects has started to change. Earlier developments were mostly short-duration systems, typically able to discharge at full power for around one hour. However, newer projects are increasingly designed to run for longer. In 2024, around 67\% of newly installed capacity had a two-hour duration, showing a clear shift towards assets that can provide energy over longer periods, not just fast response services~\cite{modo2025yearreview}. This trend has continued into 2025. By the end of the year, total installed battery capacity had reached around 6.8 GW, making it the strongest year of growth on record for the sector.

\subsection{Looking Ahead to 2050}
Battery storage currently makes up approximately 8\% share of Great Britain's electricity capacity, but it is set to grow rapidly. The National Energy System Operator's Future Energy Scenarios (FES) 2025 set out credible pathways to net zero by 2050~\cite{neso2025fes}. Across all pathways, storage becomes a central part of the system. Grid connected battery capacity alone is expected to reach between 31 and 40 GW by 2050. At the same time, long duration storage, such as pumped hydro and emerging technologies, is projected to grow to between 13 and 17 GW. In one of the FES' scenario (Holistic Transition), the total storage capacity, including vehicle to grid, reaches around 96 GW by 2050. This highlights the scale of flexibility required. Meeting nearer term targets is also demanding. Delivering Clean Power by 2030 requires around 3 GW of new battery capacity each year, about twice the highest annual rate seen so far.

Electric vehicles are expected to play a defining role. By 2050, there could be up to 37.4 million electric vehicles on the road. Assuming an average battery size of 70 kilowatt hours, this would amount to more than two and a half terawatt hours of installed storage. Only a share of this will be available at any time, but even limited participation through smart charging and vehicle to grid could have a major impact.


Taken together, these trends point to a fundamental shift. Grid connected batteries will provide reliable services from central locations, supporting frequency response and energy balancing. However, electric vehicle batteries, due to their scale and distribution, offer a much larger and more flexible resource. They can absorb surplus renewable generation and return it when demand is high, while also supporting local networks.


\section{Market Mechanisms and Revenue Streams}
A key lesson from the British experience is that technology alone does not deliver deployment. It is the market design that makes it viable at scale. Batteries do not generate electricity. They create value by shifting energy in time and by providing services the system needs. Their success therefore depends on whether markets recognise and reward these services.

In the early years, most revenue came from fast frequency response. Batteries could react within seconds, making them well suited to stabilising system frequency. However, as more projects connected, these markets quickly became saturated. Revenues declined and operators had to adapt. Today, battery projects rely on a stack of income sources. These include the Balancing Mechanism, wholesale trading, and the Capacity Market, which rewards availability during periods of system stress. This shift reflects a more mature sector, where value comes from combining services rather than relying on a single product.

The role of market design in driving this transition has been decisive. The introduction of Dynamic Containment in 2020 is a clear example~\cite{CAO2024110288}. Strong early prices attracted rapid investment and helped establish batteries as a core part of system operation. Subsequent services, such as Dynamic Moderation and Dynamic Regulation, expanded the range of opportunities and allowed better matching of battery capabilities to system needs~\cite{pacificgreen2024capacity}.

The Capacity Market has also been important. It has provided longer term revenue visibility and, through flexible rules, allowed battery projects to compete alongside more traditional assets. At the same time, increased participation in the Balancing Mechanism and energy markets has opened up new sources of value~\cite{modo2025yearreview}.

\section{The Emerging Role of Electric Vehicles}
\begin{figure}
\begin{tikzpicture}
\begin{axis}[
    width=\columnwidth,
    xlabel={Year},
    ylabel={Number of EVs (millions)},
    xmin=2025, xmax=2050,
    ymin=-3500, ymax=41000000,
    xtick={2025,2030,2035,2040,2045,2050},
    xticklabels={2025,2030,2035,2040,2045,2050},
    ytick={5000000,10000000,15000000,20000000,25000000,30000000,35000000,40000000},
    yticklabels={5,10,15,20,25,30,35,40},scaled ticks=false,
    legend style={
        at={(0.425,0.325)},
        anchor=north west,
        font=\small,
    },
    grid=both,
    grid style={line width=0.3pt, draw=gray!30},
    major grid style={line width=0.5pt, draw=gray!50},
    minor tick num=1,
    tick align=outside,
    axis line style={-},
    every axis plot/.append style={line width=1.6pt, mark size=1.8pt},
]
 
\addplot[color=blue, mark=circle, solid]
coordinates {
    (2025, 2005061)
    (2026, 2812035)
    (2027, 3799613)
    (2028, 4995793)
    (2029, 6426931)
    (2030, 8212830)
    (2031, 10263759)
    (2032, 12545455)
    (2033, 14999777)
    (2034, 17676549)
    (2035, 20346070)
    (2036, 22923657)
    (2037, 25291707)
    (2038, 27383207)
    (2039, 29190658)
    (2040, 30712427)
    (2041, 31938369)
    (2042, 32927098)
    (2043, 33219164)
    (2044, 33303012)
    (2045, 33333028)
    (2046, 33277382)
    (2047, 32958495)
    (2048, 32551494)
    (2049, 32080227)
    (2050, 31519388)
};
\addlegendentry{Holistic Transition}
 
\addplot[color=green, mark=square, dashed]
coordinates {
    (2025, 2005061)
    (2026, 2812036)
    (2027, 3799615)
    (2028, 4995804)
    (2029, 6426974)
    (2030, 8212963)
    (2031, 10264126)
    (2032, 12546405)
    (2033, 15002083)
    (2034, 17674952)
    (2035, 20344533)
    (2036, 22924980)
    (2037, 25300362)
    (2038, 27406378)
    (2039, 29237059)
    (2040, 30792042)
    (2041, 32061453)
    (2042, 33103427)
    (2043, 33789743)
    (2044, 33980232)
    (2045, 34138295)
    (2046, 34259743)
    (2047, 34122695)
    (2048, 33931028)
    (2049, 33694799)
    (2050, 33399709)
};
\addlegendentry{Electric Engagement}
 
\addplot[color=purple, mark=triangle, dashdotted]
coordinates {
    (2025, 2005368)
    (2026, 2805977)
    (2027, 3784954)
    (2028, 4968705)
    (2029, 6384383)
    (2030, 8161759)
    (2031, 10204977)
    (2032, 12479945)
    (2033, 14929824)
    (2034, 17599894)
    (2035, 20273439)
    (2036, 22850662)
    (2037, 25224746)
    (2038, 27339217)
    (2039, 29213364)
    (2040, 30804671)
    (2041, 32156088)
    (2042, 33293057)
    (2043, 34251936)
    (2044, 34832755)
    (2045, 35170478)
    (2046, 35514007)
    (2047, 35668357)
    (2048, 35796942)
    (2049, 35924418)
    (2050, 36050582)
};
\addlegendentry{Hydrogen Evolution}
 
\addplot[color=red, mark=diamond, dotted]
coordinates {
    (2025, 1768392)
    (2026, 2356394)
    (2027, 3042265)
    (2028, 3841829)
    (2029, 4770415)
    (2030, 5827136)
    (2031, 7023768)
    (2032, 8336102)
    (2033, 9758078)
    (2034, 11261656)
    (2035, 12875742)
    (2036, 14747383)
    (2037, 16675235)
    (2038, 18609833)
    (2039, 20519512)
    (2040, 22648712)
    (2041, 24717891)
    (2042, 26655492)
    (2043, 28480186)
    (2044, 30176099)
    (2045, 31732223)
    (2046, 32849839)
    (2047, 33422843)
    (2048, 34026289)
    (2049, 34634706)
    (2050, 35174313)
};
\addlegendentry{Falling Behind}
 
\end{axis}
\end{tikzpicture}
\caption{Projected number of electric vehicles in Great Britain under the four FES scenarios. Data from~\cite{neso2025fes}}
\label{fig:fesevs}
\end{figure}
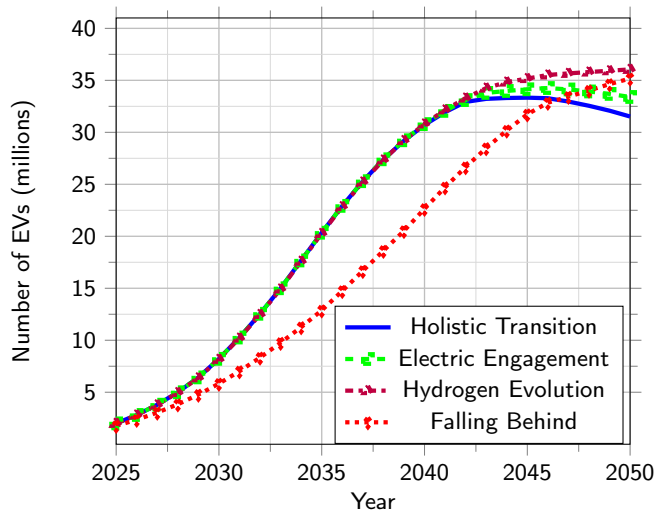

The NESO's Future Energy Scenarios estimate that between 30 and 36 million electric vehicles will be on Great Britain's roads by 2050 to meet net zero targets. Figure~\ref{fig:fesevs} presents the evolution of electric vehicles under the four scenarios up to the year 2050. With typical battery sizes between 50 and 100 kilowatt hours, even a conservative assumption of 31 million vehicles with 60 kilowatt hours each would give a total of around 1,860 gigawatt hours of storage.

The level of flexibility may be offered by electric vehicles is orders of magnitude larger than today's grid connected battery capacity, which is only a few tens of gigawatt hours. Within the lifetime of vehicles already being sold, the largest source of energy storage in Great Britain will no longer be on the transmission system, but in homes, workplaces, and depots. This shift is reflected in system planning. The National Energy System Operator identifies electric vehicles as the single largest source of flexibility by 2050. In its analysis, electric vehicles could provide up to 51 GW of flexible capacity, exceeding the capacity of today's gas fired generation fleet.

Vehicle to grid technology is central to this role. It allows electricity to flow both to and from vehicles, turning them into active system resources. Updated scenarios now include not only private cars but also commercial fleets such as vans, buses, and heavy goods vehicles. This increases the estimated vehicle to grid capacity to around 41 GW in the Holistic Transition pathway. Importantly, this flexibility can be provided without affecting users, as vehicles can still be fully charged when needed.

\subsection{The Coordination Challenge}
The potential flexibility offered by electric vehicles is very large, but it is important to be clear about when this becomes a benefit to the system and when it becomes a burden. If electric vehicles simply start charging as soon as they are plugged in, they tend to increase demand at the worst possible time. This is usually in the early evening, when the system is already under pressure. In this situation, they are not providing flexibility but instead adding stress.

Studies consistently show that unmanaged charging increases peak demand and can lead to local voltage problems and network overload~\cite{abdelhamid2024evimpacts}. The difference between a problem and a solution lies in control. Smart charging, which shifts demand to times when the system is under less pressure, is already becoming standard practice. It helps to reduce peaks and improve overall system efficiency.

However, this alone does not unlock the full value of electric vehicles. To achieve that, vehicle to grid capability is needed. When both charging and discharging are actively controlled, electric vehicles can respond to the needs of the system rather than the habits of users. At that point, the main challenge is no longer technical capability, but coordination at scale.

\begin{figure}

\begin{tikzpicture}
\begin{axis}[
    width=\columnwidth,
    xlabel={Year},
    ylabel={Number of EVs (millions)},
    xmin=2025, xmax=2050,
    ymin=-1500000, ymax=3000000,
    xtick={2025,2030,2035,2040,2045,2050},
    xticklabels={2025,2030,2035,2040,2045,2050},
    ytick={-1000000,0,1000000,2000000,3000000},
    yticklabels={-1,0,1,2,3},scaled ticks=false,
    legend style={
        at={(0.1,0.4)},
        anchor=north west,
        font=\small,
    },
    grid=both,
    grid style={line width=0.3pt, draw=gray!30},
    major grid style={line width=0.5pt, draw=gray!50},
    minor tick num=1,
    tick align=outside,
    axis line style={-},
    every axis plot/.append style={line width=1.6pt, mark size=1.8pt},
]
 
\addplot[color=blue, mark=circle, solid]
coordinates {
    (2026,  806974)
    (2027,  987578)
    (2028,  1196180)
    (2029,  1431138)
    (2030,  1785899)
    (2031,  2050929)
    (2032,  2281696)
    (2033,  2454322)
    (2034,  2676772)
    (2035,  2669521)
    (2036,  2577587)
    (2037,  2368050)
    (2038,  2091500)
    (2039,  1807451)
    (2040,  1521769)
    (2041,  1225942)
    (2042,  988729)
    (2043,  292066)
    (2044,  83848)
    (2045,  30016)
    (2046,  -55646)
    (2047,  -318887)
    (2048,  -407001)
    (2049,  -471267)
    (2050,  -560839)
};
\addlegendentry{Holistic Transition}
 
\addplot[color=green, mark=square, dashed]
coordinates {
    (2026,  806975)
    (2027,  987579)
    (2028,  1196189)
    (2029,  1431170)
    (2030,  1785989)
    (2031,  2051163)
    (2032,  2282279)
    (2033,  2455678)
    (2034,  2672869)
    (2035,  2669581)
    (2036,  2580447)
    (2037,  2375382)
    (2038,  2106016)
    (2039,  1830681)
    (2040,  1554983)
    (2041,  1269411)
    (2042,  1041974)
    (2043,  686316)
    (2044,  190489)
    (2045,  158063)
    (2046,  121448)
    (2047,  -137048)
    (2048,  -191667)
    (2049,  -236229)
    (2050,  -295090)
};
\addlegendentry{Electric Engagement}
 
\addplot[color=purple, mark=triangle, dashdotted]
coordinates {
    (2026,  800609)
    (2027,  978977)
    (2028,  1183751)
    (2029,  1415678)
    (2030,  1777376)
    (2031,  2043218)
    (2032,  2274968)
    (2033,  2449879)
    (2034,  2670070)
    (2035,  2673545)
    (2036,  2577223)
    (2037,  2374084)
    (2038,  2114471)
    (2039,  1874147)
    (2040,  1591307)
    (2041,  1351417)
    (2042,  1136969)
    (2043,  958879)
    (2044,  580819)
    (2045,  337723)
    (2046,  343529)
    (2047,  154350)
    (2048,  128585)
    (2049,  127476)
    (2050,  126164)
};
\addlegendentry{Hydrogen Evolution}
 
\addplot[color=red, mark=diamond, dotted]
coordinates {
    (2026,  588002)
    (2027,  685871)
    (2028,  799564)
    (2029,  928586)
    (2030,  1056721)
    (2031,  1196632)
    (2032,  1312334)
    (2033,  1421976)
    (2034,  1503578)
    (2035,  1614086)
    (2036,  1871641)
    (2037,  1927852)
    (2038,  1934598)
    (2039,  1909679)
    (2040,  2129200)
    (2041,  2069179)
    (2042,  1937601)
    (2043,  1824694)
    (2044,  1695913)
    (2045,  1556124)
    (2046,  1117616)
    (2047,  573004)
    (2048,  603446)
    (2049,  608417)
    (2050,  539607)
};
\addlegendentry{Falling Behind}
 
\end{axis}
\end{tikzpicture}
\caption{Annual rate of change in the number of electric vehicles to 2050 under the four FES scenarios.}
\label{fig:diffEV}
\end{figure}
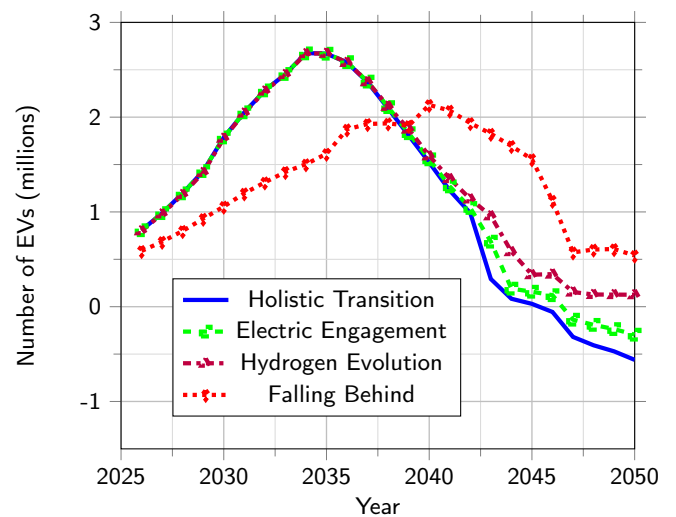

Figure~\ref{fig:diffEV} shows the rate of increase in electric vehicles under the four FES scenarios. The number of new vehicles rises each year and reaches a peak around 2034 in the more ambitious scenarios, before starting to decline. Adding around one million electric vehicles each year until 2034 means that a whole system approach is needed to integrate them effectively. This includes charging infrastructure, upgrades to distribution networks, and additional capacity to support charging demand. However, if charging is well coordinated, the need for major system upgrades can be reduced.

\subsection{The Aggregation Problem}
The scale and structure of the electric vehicle fleet make coordination fundamentally different from traditional storage. A grid connected battery is a single, controllable asset. Its location, capacity, and response are known and predictable. An electric vehicle fleet is the opposite. It consists of millions of small devices, each connected at the distribution level, each with its own usage pattern, and each appearing and disappearing from the system at different times.

This creates two problems. First, system level coordination becomes complex. Second, poorly coordinated behaviour can shift problems from the transmission system to local networks. A fully centralised approach, where a single operator controls every vehicle in real time, is not realistic. The communication burden would be excessive, and delays would undermine performance. Instead, the emerging solution is layered control.

Fast decisions are made locally at the level of the vehicle or charger. These are coordinated at site or network level to manage local constraints. At a higher level, aggregators combine many vehicles into a single market facing resource and respond to system signals over longer timescales. Aggregators are therefore a critical part of the future system. They translate system needs into large numbers of small actions. As electric vehicles scale, this layer becomes essential.

Unlocking this capability depends on advances in optimisation, control, and data driven decision making. Efficient scheduling, participation across multiple markets, and managing uncertainty in demand and renewable output are all active areas of research and deployment. Without these tools, the theoretical flexibility of electric vehicles cannot be realised in practice.

\subsection{Battery Degradation and User Trust}
Technical coordination alone is not enough. The success of vehicle to grid also depends on user participation. A key concern is battery degradation. Additional charging and discharging cycles can increase wear, and users are sensitive to any impact on battery life and warranties. This concern is justified, but it is often overstated.

Evidence shows that degradation depends more on how batteries are used than how often they are used. Deep discharge and poor thermal conditions drive most of the damage. In contrast, shallow and well managed cycling, which is what coordinated vehicle to grid would require, has a much smaller impact. This means that degradation can be limited through careful control. However, this must be communicated clearly and backed by credible guarantees. Without trust, participation will remain limited.

Incentives are equally important. Users must see clear and consistent value from participation. Early programmes have shown that this is achievable. When financial rewards are visible and simple, households and businesses are willing to engage. As electric vehicles become central to system operation, market design will need to reflect this. Flexibility services must be accessible to distributed resources, and revenues must flow back to users in a transparent way. Without this, the largest storage resource in the system will remain underused.

\section{Conclusion}
Batteries are becoming an ever larger and more important part of the Great Britain energy system. From a standing start little more than a decade ago, grid connected storage has already reached multi gigawatt scale and is set to grow several times over. Electric vehicle batteries will take this growth to another level entirely. Their combined capacity will far exceed dedicated storage plants and provide the flexible backbone needed for a renewable led power system.

The biggest challenge ahead is to streamline the participation of electric vehicles. Each battery is small and spread across a wide geographical area, owned by millions of individuals and businesses with different driving patterns and priorities. Coordinating them is not simple. If we leave charging unchecked, electric vehicles could create new evening peaks that strain the grid and increase costs. But with smart algorithms, real time pricing signals, and vehicle to grid technology, this challenge is not impossible. It is achievable. Coordinated dispatch can turn millions of small batteries into a single, powerful resource that delivers system value far beyond what centralised storage alone can provide.



\begin{thebibliography}{18}
\providecommand{\natexlab}[1]{#1}
\providecommand{\url}[1]{\texttt{#1}}
\expandafter\ifx\csname urlstyle\endcsname\relax
  \providecommand{\doi}[1]{doi: #1}\else
  \providecommand{\doi}{doi: \begingroup \urlstyle{rm}\Url}\fi

\bibitem[{United Nations Framework Convention on Climate
  Change}(2015)]{unfccc2015paris}
{United Nations Framework Convention on Climate Change}.
\newblock The {Paris} agreement, 2015.
\newblock URL
  \url{https://unfccc.int/process-and-meetings/the-paris-agreement}.

\bibitem[{International Renewable Energy Agency}(2024)]{irena2024renewable}
{International Renewable Energy Agency}.
\newblock Renewable capacity statistics 2024.
\newblock Technical report, IRENA, Abu Dhabi, 2024.
\newblock URL
  \url{https://www.irena.org/Publications/2024/Mar/Renewable-capacity-statistics-2024}.

\bibitem[Chu and Majumdar(2012)]{chu2012opportunities}
Steven Chu and Arun Majumdar.
\newblock Opportunities and challenges for a sustainable energy future.
\newblock \emph{Nature}, 488\penalty0 (7411):\penalty0 294--303, 2012.
\newblock \doi{10.1038/nature11475}.

\bibitem[{International Energy Agency}(2021)]{iea2021storage}
{International Energy Agency}.
\newblock Energy storage.
\newblock Technical report, IEA, Paris, 2021.
\newblock URL \url{https://www.iea.org/reports/energy-storage}.

\bibitem[Guerra et~al.(2024)Guerra, Welfle, Gutiérrez-Alvarez, Freer, Ma, and
  Haro]{GUERRA2024122447}
K.~Guerra, A.~Welfle, R.~Gutiérrez-Alvarez, M.~Freer, L.~Ma, and P.~Haro.
\newblock The role of energy storage in great britain's future power system:
  focus on hydrogen and biomass.
\newblock \emph{Applied Energy}, 357:\penalty0 122447, 2024.
\newblock ISSN 0306-2619.
\newblock \doi{https://doi.org/10.1016/j.apenergy.2023.122447}.
\newblock URL
  \url{https://www.sciencedirect.com/science/article/pii/S0306261923018111}.

\bibitem[{BloombergNEF}(2023)]{bnef2023evo}
{BloombergNEF}.
\newblock Electric vehicle outlook 2023.
\newblock Technical report, BloombergNEF, New York, 2023.
\newblock URL \url{https://about.bnef.com/electric-vehicle-outlook/}.

\bibitem[{The Royal Society}(2023)]{RoyalSociety2023LargeScale}
{The Royal Society}.
\newblock Large-scale electricity storage.
\newblock Report, The Royal Society, London, September 2023.
\newblock URL
  \url{https://royalsociety.org/-/media/policy/projects/large-scale-electricity-storage/large-scale-electricity-storage-report.pdf}.
\newblock DES6851\_1.

\bibitem[{International Energy Agency}(2024)]{iea2024battery}
{International Energy Agency}.
\newblock Batteries and secure energy transitions.
\newblock Technical report, IEA, Paris, 2024.
\newblock URL
  \url{https://www.iea.org/reports/batteries-and-secure-energy-transitions}.

\bibitem[{Department for Energy Security and Net
  Zero}(2024)]{DESNZ2024_clean_power_2030}
{Department for Energy Security and Net Zero}.
\newblock Clean power 2030 action plan: A new era of clean electricity.
\newblock Technical report, UK Government, December 2024.
\newblock URL
  \url{https://www.gov.uk/government/publications/clean-power-2030-action-plan}.
\newblock Published 13 December 2024; accessed 2026-04-09.

\bibitem[{National Energy System Operator
  (NESO)}(2025)]{NESO_AnnualBalancingReport_2025}
{National Energy System Operator (NESO)}.
\newblock Annual balancing report 2024/25.
\newblock Technical report, National Energy System Operator, June 2025.
\newblock URL \url{https://www.neso.energy/document/362561/download}.
\newblock Accessed: 2026-04-23.

\bibitem[{National Energy System Operator}(2025)]{neso2025fes}
{National Energy System Operator}.
\newblock Future energy scenarios: Pathways to net zero 2025.
\newblock Technical report, National Energy System Operator ({NESO}), July
  2025.
\newblock URL
  \url{https://www.neso.energy/publications/future-energy-scenarios-fes}.
\newblock Accessed: April 2026.

\bibitem[{UK Power Networks}(2016)]{UKPN_SNS_2016}
{UK Power Networks}.
\newblock Smarter network storage (sns).
\newblock
  \url{https://innovation.ukpowernetworks.co.uk/projects/smarter-network-storage-sns},
  2016.
\newblock Accessed: 2026-04-27.

\bibitem[{National Grid ESO}(2016)]{nationalgrid2016efr}
{National Grid ESO}.
\newblock Enhanced frequency response: Frequently asked questions.
\newblock Technical report, National Grid Electricity System Operator, 2016.
\newblock URL
  \url{https://www.nationalgrid.com/sites/default/files/documents/Enhanced%20Frequency%20Response%20FAQs%20v5.0_.pdf}.
\newblock Accessed: April 2026.

\bibitem[Stoker(2019)]{stoker2019blackout}
Liam Stoker.
\newblock Long read: How energy storage helped save the {UK}'s grid.
\newblock \emph{Solar Power Portal}, August 2019.
\newblock URL
  \url{https://www.solarpowerportal.co.uk/residential-solar/long-read-how-energy-storage-helped-save-the-uk-s-grid}.
\newblock Accessed: April 2026.

\bibitem[{Modo Energy}(2025)]{modo2025yearreview}
{Modo Energy}.
\newblock {GB} battery energy storage markets: 2024 year in review.
\newblock Technical report, Modo Energy, January 2025.
\newblock URL \url{https://tinyurl.com/bdh62j9u}.
\newblock Accessed: April 2026.

\bibitem[Cao et~al.(2024)Cao, Engelhardt, Ziras, Marinelli, and
  Zhao]{CAO2024110288}
Xihai Cao, Jan Engelhardt, Charalampos Ziras, Mattia Marinelli, and Nan Zhao.
\newblock Battery energy storage systems providing dynamic containment
  frequency response service.
\newblock \emph{International Journal of Electrical Power \& Energy Systems},
  162:\penalty0 110288, 2024.
\newblock ISSN 0142-0615.
\newblock \doi{https://doi.org/10.1016/j.ijepes.2024.110288}.
\newblock URL
  \url{https://www.sciencedirect.com/science/article/pii/S0142061524005106}.

\bibitem[{Pacific Green}(2024)]{pacificgreen2024capacity}
{Pacific Green}.
\newblock Is capacity the future for {UK} battery storage?
\newblock Pacific Green, 2024.
\newblock URL
  \url{https://www.pacificgreen.com/articles/capacity-future-uk-battery-storage/}.
\newblock Accessed: April 2026.

\bibitem[Abdelhamid et~al.(2024)]{abdelhamid2024evimpacts}
Mohamed Abdelhamid et~al.
\newblock Analysis of multidimensional impacts of electric vehicles penetration
  in distribution networks.
\newblock \emph{Scientific Reports}, 14:\penalty0 27854, November 2024.
\newblock \doi{10.1038/s41598-024-77662-6}.
\newblock URL \url{https://www.nature.com/articles/s41598-024-77662-6}.

\end{thebibliography}

\end{document}